\documentclass[review]{elsarticle}

\usepackage{lineno,hyperref}
\modulolinenumbers[5]
\usepackage{multirow}
\usepackage{makerobust}
\MakeRobustCommand\overrightarrow
\MakeRobustCommand\overleftarrow
\usepackage{pdflscape}
\usepackage{endnotes}
\usepackage{amsmath}
\usepackage{graphicx}
\usepackage{caption}
\usepackage{subcaption}
\usepackage{amsmath, amsthm, amssymb}
\usepackage{appendix}
\usepackage{graphicx}
\usepackage{resizegather}
\usepackage{color}
\usepackage{xcolor}
\newtheorem{proposition}{Proposition}
\newtheorem{corollary}{Corollary}
\newtheorem{lemma}{Lemma}

\begin{document}

\begin{frontmatter}

\title{Approximate Performance Measures for a Single Station Two-Stage Reneging Queue}


\author[mymainaddress]{Katsunobu Sasanuma\corref{mycorrespondingauthor}}
\ead{katsunobu.sasanuma@stonybrook.edu}

\author[mysecondaryaddress]{Alan Scheller-Wolf}
\cortext[mycorrespondingauthor]{Corresponding author}
\ead{awolf@andrew.cmu.edu}

\address[mymainaddress]{College of Business, Stony Brook University, Stony Brook, NY 11794, USA}
\address[mysecondaryaddress]{Tepper School of Business, Carnegie Mellon University, Pittsburgh, PA 15213, USA}

\begin{abstract}
We study a single station \emph{two-stage} reneging queue with Poisson arrivals, exponential services, and two levels of exponential reneging behaviors, extending the popular Erlang A model that assumes a constant reneging rate. We derive approximate analytical formulas representing performance measures for our two-stage queue following the Markov chain decomposition approach. Our formulas not only give accurate results spanning the heavy-traffic to the light-traffic regimes, but also provide insight into capacity decisions.
\end{abstract}


\end{frontmatter}

\linenumbers

\section{Introduction}
Models representing service systems with impatient customers have been studied extensively due to their practical importance. When analyzing such systems, it is often assumed that all customers are equally impatient and randomly renege (i.e., leave a queue after entering but before reaching service) at the same rate. This constant patience assumption may be reasonable when customers have no information about their position in queue (e.g., in some call centers). However, customers often change their level of patience, for example based on their queue position if they experience a different level of comfort, or when they know they are close to, or far from, the front of the queue \cite{kuzu2019wait}. For example, many restaurants provide some seats (capacity) for customers waiting inside (1\textsuperscript{st} stage queue), but not everybody can sit there; if all seats inside are filled, any remaining waiting customers must wait outside where it may be cold or raining (2\textsuperscript{nd} stage queue). Our model aims to provide simple, yet accurate formulas for practitioners who make capacity (in this case waiting room seating capacity) decisions when their customers show stage-dependent reneging behaviors. Specifically, we consider a single station two-stage reneging queue with Poisson arrivals, exponential services, and stage-dependent exponential reneging rates. We analyze three performance measures: probability of queueing, probability of customer abandonment (via reneging or blocking), and average queue length. Our two-stage reneging model can cope with both finite and infinite queues.

When the two reneging rates match, our queue reduces to the Erlang A model, an M/M/n+M queueing model with reneging. This Erlang A model and its variations have been utilized for the analysis of various real-world problems such as public housing \cite{kaplan1987analyzing}, kidney transplantation \cite{zenios1999modeling}, on-street parking \cite{larson2010congestion}, and call centers \cite[see, e.g.,][]{garnett2002designing,gans2003telephone,borst2004dimensioning, mandelbaum2007service,baron2009staffing}. Most of the studies are based on either heavy traffic approximation or asymptotic analysis \cite{halfin1981heavy}, which are very effective for analyzing a single stage reneging queue. For more general queueing systems with state-dependent general service times, the queueing and Markov chain decomposition (QMCD) approach is versatile \cite[see, e.g.,][]{abouee2016state,baron2018state}; the QMCD method can be used to analyze complicated models, such as a queue with orbit or non-Markovian systems. Since our model has only two stages, we utilize a simpler Markov chain decomposition approach \cite{sasanuma2019markov}, and evaluate each decomposed sub-system separately using a Poisson-Normal approximation \cite{tijms2003first} with a continuity correction term \cite{feller1968probability,cox1970continuity,tijms1986stochastic}. This continuity correction term makes our formulas accurate and robust over a wide range of parameters.

\section{Single Station Two-Stage Reneging Model}
\label{sec:two-stage}
In this section, we present our model and derive several basic performance measures following the Markov chain decomposition approach \cite{sasanuma2019markov}. The formulas we derive in this section are exact; approximations are discussed in the following section.

\subsection{Markov Chain Structure and Its Decomposed Subchains}
\label{sec:model}
\begin{figure}
	\caption{Two-Stage Reneging Model: Full Markov Chain.}
	\label{fig:fullMC}
	\centering
	\includegraphics*[scale=0.42]{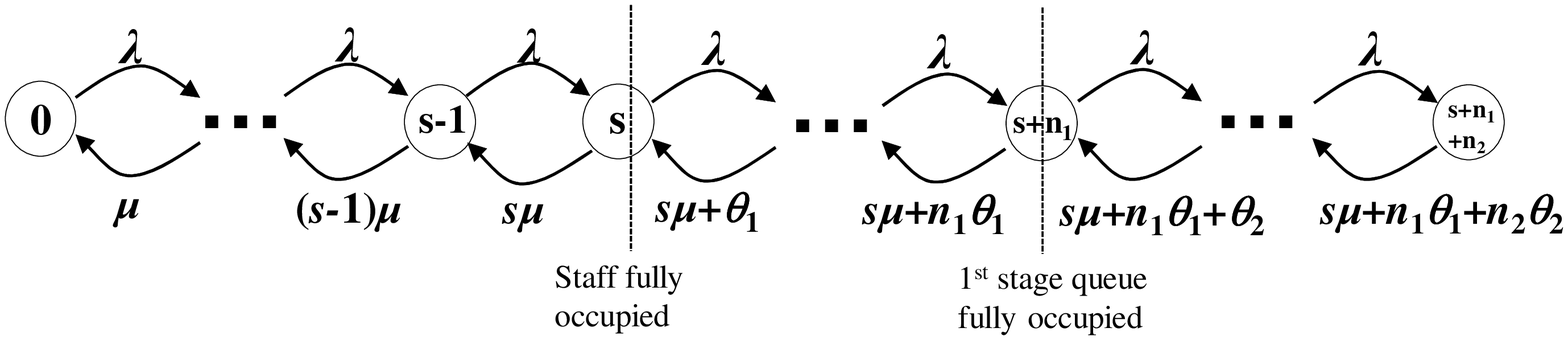}
\end{figure}

\begin{figure}
    \caption{Decomposed Markov chains.}
        \label{fig:SC}
    \centering
    \begin{subfigure}[]{0.40\textwidth}
     \centering
        \includegraphics[scale=0.4]{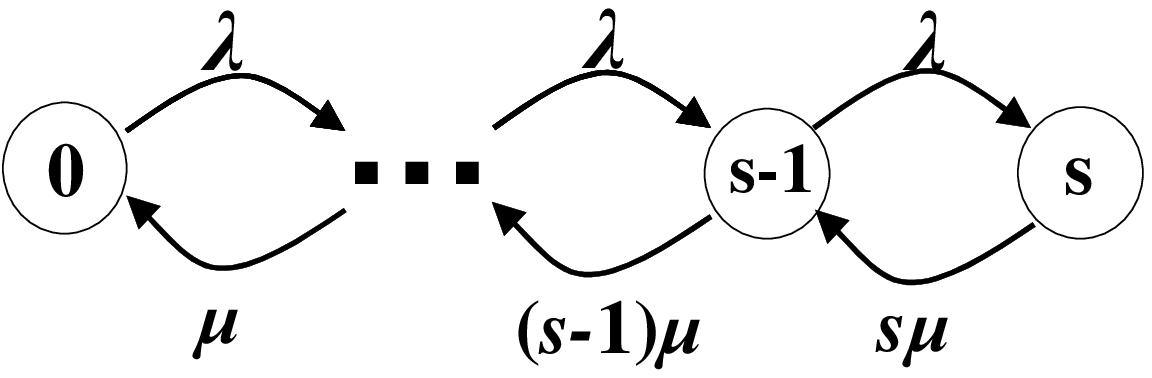}
        \caption{Subchain 0}
        \label{fig:SC0}
   \end{subfigure}
   \begin{subfigure}[]{0.58\textwidth}
     \centering
        \includegraphics[scale=0.4]{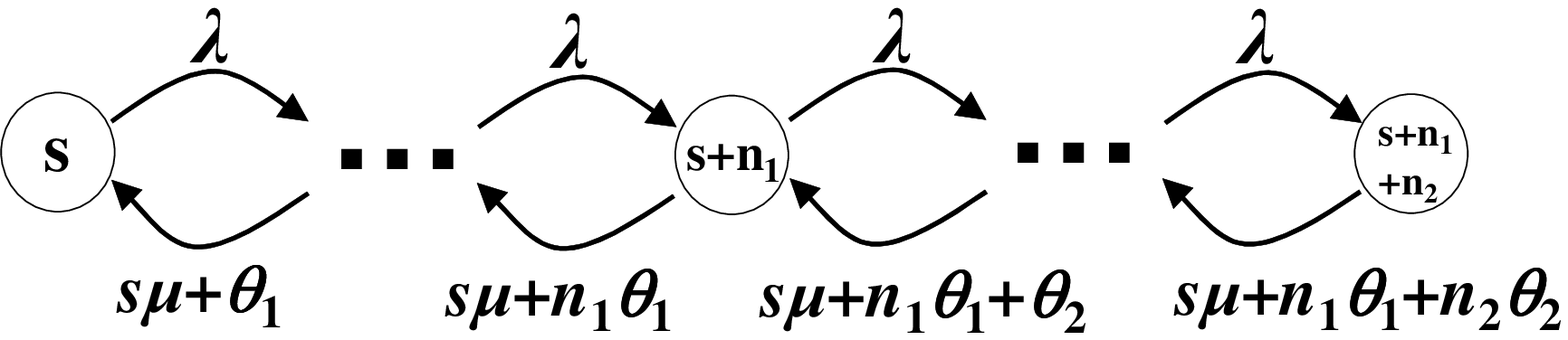}
        \caption{Subchain $q$}
        \label{fig:SCq}
    \end{subfigure}
     \begin{subfigure}[]{0.31\textwidth}
     \centering
        \includegraphics[scale=0.42]{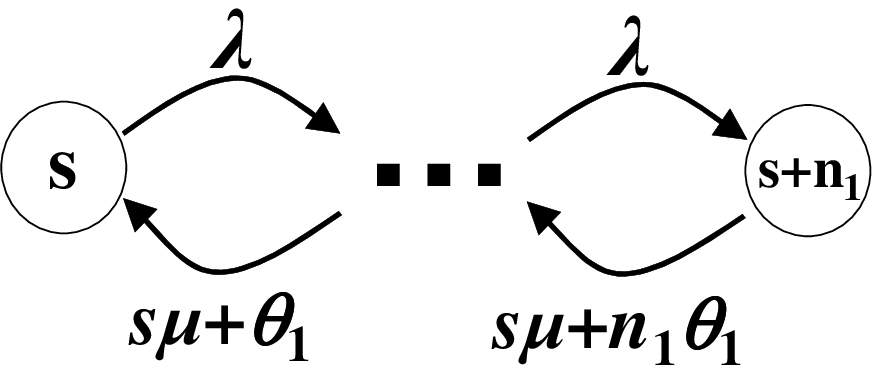}
        \caption{Subchain 1}
        \label{fig:SC1}
    \end{subfigure}
    \begin{subfigure}[]{0.38\textwidth}
     \centering
        \includegraphics[scale=0.42]{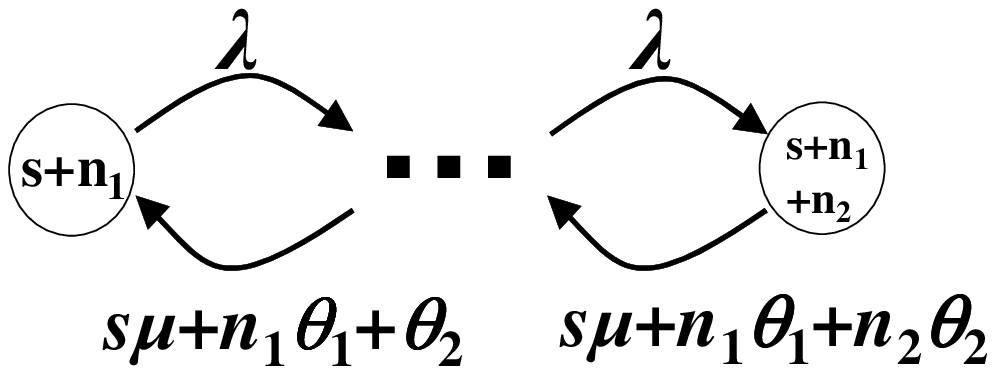}
        \caption{Subchain 2}
        \label{fig:SC2}
    \end{subfigure}
\end{figure}

Our two-stage model is a simple birth-death Markov chain (MC) (Fig.~\ref{fig:fullMC}). Customers arrive according to Poisson process with rate $\lambda$, and are served by one of the $s$ servers with exponential rate $\mu$ following a first-in first-out (FIFO) rule. This single station queue has two stages ($i=1,2$); stage $i$ has $n_i$ spaces; $n_i$ can be either finite or infinite. Each customer in the $i$\textsuperscript{th} stage reneges after an exponentially distributed time with rate $\theta_i \; (>0)$. When both stages of the queue are fully occupied in state $s+n_1+n_2$, a new customer is unable to enter the system (\emph{blocked}). The arrival and departure rates at state $k$ are
\small
\begin{eqnarray*}
\lambda _{k}&=&\left\{\begin{array}{ll} {\lambda } {\qquad}{\qquad}& {0 \le k<s+n_1+n_2} \\ {0} & {s+n_1+n_2\le k} \end{array}\right. \\
\mu _{k} &=&\left\{\begin{array}{ll} {k\mu } & {1 \le k\le s} \\ {s\mu +\left(k-s\right)\theta_1} & {s<k\le s+n_1} \\ {s\mu + n_1\theta_1+\left(k-s-n_1 \right)\theta_2}& {s+n_1<k\le s+n_1+n_2}\\ {0} & {s+n_1+n_2 < k.} \end{array}\right.
\end{eqnarray*}
\normalsize
To analyze our system we decompose the full MC (Fig.~\ref{fig:fullMC}) into subchains: an $M/M/s/s$ subchain (Subchain 0 in Fig.~\ref{fig:SC0}: MC comprised of states 1 to $s$) and a queueing subchain (Subchain $q$ in Fig.~\ref{fig:SCq}:MC comprised of states $s$ to $s+n_1+n_2$). Subchain $q$ is further decomposed into the 1\textsuperscript{st} stage queue with reneging rate $\theta_1$ (Subchain 1 in Fig.~\ref{fig:SC1}: MC comprised of states $s$ to $s+n_1$) and the 2\textsuperscript{nd} stage queue with reneging rate $\theta_2$ (Subchain $2$ in Fig.~\ref{fig:SC2}: MC comprised of states $s+n_1$ to $s+n_1+n_2$). Let $\Omega$ and $A_j$ represent a set of states for the full MC and subchain $j(=0,1,2,q)$, respectively. Thus: $\Omega=A_0 \cup A_q$, $A_q=A_1 \cup A_2$, $A_0 \cap A_q = A_0 \cap A_1 = \{s\}$, and $A_1 \cap A_2 = \{s+n_1\}$.

\subsection{Representation of Performance Measures}
\label{sec:representation}
We denote the steady-state probability of state $k$ in the full MC and subchain $j$ as $\pi _{k}$ and $\pi_k^j$, respectively. Since the MC we consider is of the birth-death type, any truncated subchain $j$ maintains the same stationary distribution as the full MC up to a normalization constant, i.e., $\pi_k^j=\frac{\pi_k}{\sum_{k' \in A_j} \pi_{k'}}, \forall k \in A_j$ (which we simply denote as $\pi_k^j \propto \pi_k, \forall k \in A_j$ in this paper). Define $\mathbb{E}[{f(X)}]=\sum_{k \in \Omega} f(k)\pi_k$, $\mathbb{E}_j[{f(X)}]=\sum_{k \in A_j} f(k)\pi^j_k$, and ${r_1 = \frac{\pi^1_{s+n_1}}{\pi^1_{s}}}$. The following proposition holds as a special case of the Markov chain decomposition method based on the total expectation theorem \cite{sasanuma2019markov}. All proofs are given in appendices.
\begin{proposition}
\label{tpt}
For any function of states $f(X)$ for the two-stage reneging queue, the following equations hold:
\small
\begin{equation}
\label{tpt1}
\frac{\mathbb{E}[ {f\left( X \right)}]}{\pi_s}=\frac{\mathbb{E}_0[ {f\left( X \right)}]}{\pi^0_s}+\frac{\mathbb{E}_q[ {f\left( X \right)}]}{\pi^q_s}-f(s),\\
\label{tpt2}
\frac{\mathbb{E}_q[ {f\left( X \right)}]}{\pi^q_s}=\frac{\mathbb{E}_1[ {f\left( X \right)}]}{\pi^1_s}+r_1\left(\frac{\mathbb{E}_2[ {f\left( X \right)}]}{\pi^2_{s+n_1}}-f(s+n_1)\right).
\end{equation}
\normalsize
Combining above, we obtain
\small
\begin{equation}
\label{tpt3}
\frac{\mathbb{E}[ {f\left( X \right)}]}{\pi_s}=\frac{\mathbb{E}_0[ {f\left( X \right)}]}{\pi^0_s}+\left(\frac{\mathbb{E}_1[ {f\left( X \right)}]}{\pi^1_s}-f(s)\right)+r_1 \left(\frac{\mathbb{E}_2[ {f\left( X \right)}]}{\pi^2_{s+n_1}}-f(s+n_1)\right).
\end{equation}
\normalsize
\end{proposition}

Denote the queueing probability as $P_Q$, the abandonment probability (due to reneging and blocking) as $P_A$, and the expected number of customers in queue as $L$. Using Proposition~\ref{tpt}, we can obtain $1/\pi_s$ by setting $f(X)=1$; $P_Q$ by setting $f(X)=\bold{1}_{A_q}$ (an indicator function); $P_A$ by setting $f(X)=N_A$ (a random variable representing the steady-state number of customer abandonments per unit time); and $L$ by setting $f(X)=N$ (a random variable representing the steady-state number of customers in queue).
\begin{proposition}
\label{str-rep}
The performance measures of the two-stage reneging queue are
\small
\begin{align}
\label{pi}
\frac{1}{\pi_s} &=\frac{1}{\pi _s^0}+\frac{1}{\pi _s^q}-1=\frac{1}{\pi _s^0}+ \left( \frac{1}{\pi _s^1}-1 \right) +r_1 \left(\frac{1}{\pi _{s+n_1}^2}-1\right),\\
\label{pq}
\frac{P_Q}{\pi_s} &= \frac{1}{\pi _s^q}=\frac{1}{\pi _s^1}+r_1\left(\frac{1}{\pi _{s+n_1}^2}-1\right),\\
\label{pa}
\frac{P_A}{\pi_s}&= p\left(\frac{1}{\pi_s^q}-1\right)+1= p\left[\left(\frac{1}{\pi _s^1}-1\right)+r_1\left(\frac{1}{\pi _{s+n_1}^2}-1\right)\right]+1,\\
\label{el}
\frac{L}{\pi_s}&=\tfrac{\lambda}{\theta_1}\left[p \left(\tfrac{1}{\pi _{s}^1}-1\right) + 1-r_1 \right]+\tfrac{r_1 \lambda}{\theta_2} \left[\left(p+\tfrac{n_1(\theta_2-\theta_1)}{\lambda}\right) \left(\tfrac{1}{\pi _{s+n_1}^2}-1\right) + 1-r_2 \right],
\end{align}
\normalsize
where
\small
\begin{equation*}
r_1 = \frac{\pi^1_{s+n_1}}{\pi^1_{s}},\;
r_2 = \frac{\pi^2_{s+n_1+n_2}}{\pi^2_{s+n_1}},\; \normalsize\text{and }\mbox{ } \small
p = 1-\frac{s\mu}{\lambda}.
\end{equation*}
\end{proposition}
\normalsize
\hspace{-0.5cm}Proposition~\ref{str-rep} expresses the performance measures of interest as simple functions of steady-state probabilities; these probabilities' exact expressions are given by
\small
\begin{equation}
\label{exact}
\frac{1}{\pi_s^0}=\frac{\sum_{k=0}^s \frac{(\lambda/\mu)^k}{k!}}{\frac{(\lambda/\mu)^s}{s!}},\; \frac{1}{\pi_s^1}=\sum_{k=0}^{n_1} \prod_{j=1}^k \frac{\lambda}{s\mu+j\theta_1},\; \frac{1}{\pi_{s+n_1}^2}=\sum_{k=0}^{n_2} \prod_{j=1}^k \frac{\lambda}{s\mu+n_1\theta_1+j\theta_2},
\end{equation}
\normalsize
where we have an empty product (when $k=0$), whose value is 1 by convention. When $n_1$ (or $n_2$)$=\infty$, the last two equations have the following alternative expressions (the derivation utilizes gamma and beta functions; see \cite{zenios1999modeling}):
\small
\begin{equation}
\label{exact2}
\frac{1}{\pi_s^1}=\frac{s\mu}{\theta_1} \int_0^1 e^{\frac{\lambda t}{\theta_1}} (1-t)^{\frac{s\mu}{\theta_1}-1} dt,\; \frac{1}{\pi_{s+n_1}^2}=\frac{s\mu + n_1\theta_1}{\theta_2} \int_0^1 e^{\frac{\lambda t}{\theta_2}} (1-t)^{\frac{s\mu + n_1\theta_1}{\theta_2}-1} dt.
\end{equation}
\normalsize
Proposition~\ref{str-rep} with \eqref{exact} and \eqref{exact2} gives the exact performance measures of our model; however, the result is cumbersome to evaluate and provides little intuition on the performance measures of interest. To gain some insight about the performance measures, the following exact relationships (directly derived from Proposition~\ref{str-rep}) are useful.
\begin{corollary}
\label{rules}
$P_A$ is related to $P_Q$ as follows:
\small
\begin{equation}
\label{pa1}
P_A = p(P_Q-\pi_s)+\pi_s.
\end{equation}
\normalsize
In addition, if $\theta=\theta_1=\theta_2$, then $L$ is related to $P_Q$ as follows:\smallskip
\small
\begin{equation}
\label{pa3}
L = \frac{\lambda}{\theta}\left[p (P_Q-\pi_s) +\pi_s - \pi_{s+n_1+n_2}\right].
\end{equation}
\end{corollary}
\normalsize
\hspace{-0.55cm}Corollary~\ref{rules} shows the impact of the discrete state space of the Markov chain on the performance measures. If we can assume that the system is sufficiently large to ignore the probabilities of being in any single state (specifically, $\pi_s$ and $\pi_{s+n_1+n_2}$ in this case), then \eqref{pa1} and \eqref{pa3} reduce to simple relationships:
$P_A \approx pP_Q \text{ and } L \approx \frac{\lambda p P_Q}{\theta},$
where $p$ is the abandonment probability when all servers are always busy. The former tells us that the probability of customer abandonment is approximately proportional to the probability of queueing. The latter is simply Little's law: the arrival to a subsystem representing a group of people who abandon the system (i.e., $\lambda p P_Q$) stays in the subsystem on average $1/\theta$ period of time.

In general, we expect that relationships among performance measures as well as analytical representations of performance measures will take simpler forms at the limit of a large system where each individual state has only a negligible contribution to the properties of the full system. Unfortunately, such a simplistic view is not appropriate when dealing with small systems, because for these systems the discreteness of the state space typically plays an important role. Thus, in our approximation procedure, we evaluate (and approximate) the discrete properties of the system without taking a limit of a large system.

\section{Approximation Procedure}
\label{sec:approximation}
Our performance measures of interest are simple functions of subchains' steady-state probabilities. Thus, to find approximate representations of the performance measures, we derive an approximate analytical representation of these steady-state probabilities.

\subsection{Steady-State Probabilities in Poisson and Normal Representations}
\label{sec:conversion}
We first show the Poisson-Normal conversion formulas, which involve a Poisson random variable (RV) with mean $R$ and the standard normal RV, whose cumulative distribution functions (CDFs) are $F_P(\cdot;R)$ and $\Phi(\cdot)$; their corresponding probability mass function (PMF) and probability density function (PDF) are $f_P(\cdot;R)$ and $\phi(\cdot)$, respectively. Define the standard normal hazard function as $h(x) = \phi(x)/(1-\Phi(x)) (=\phi(-x)/\Phi(-x))$ and introduce parameters $c={(s-R)}/{\sqrt{R}}$ and $\Delta =0.5/{\sqrt{R}}$. The following result is based on \cite{tijms1986stochastic}.
\begin{lemma}
\label{Poisson-Normal}
The Poisson CDF/PMF are approximated by the standard normal CDF/PDF as:
\small
\begin{align}
\label{cdf}
{F_P}(s;R) &\approx \Phi (c+\Delta),\\
\label{pdf}
{f_P}(s;R) &\approx \frac{\phi (c+\Delta)}{\sqrt R},\\
\label{hazard1}
\frac{f_P(s;R)}{1-F_P(s;R)} &\approx \frac{h(c +\Delta)}{\sqrt R},\\
\label{hazard2}
\frac{f_P(s;R)}{F_P(s;R)} &\approx \frac{h(-c -\Delta)}{\sqrt R}.
\end{align}
\normalsize
\end{lemma}
\smallskip
Following conventions in capacity planning, we call the argument $s$ of the Poisson PMF/CDF the \emph{staffing level}, the mean $R$ of the Poisson RV the \emph{resource requirement}, and the argument $c$ of the standard normal PDF/CDF the \emph{square-root coefficient}. The term $\Delta$ is the \emph{continuity correction term}, which is introduced due to the conversion from a discrete function (Poisson) to continuous (normal). The term $\Delta$ is non-negligible when $R$ is small (around 10 or less), but diminishes to zero in the asymptotic limit of large $R$.

\begin{table}
\begin{center}
\caption{Parameters for Poisson and Standard Normal Representations}
\small
\begin{tabular}{ cllll }
\noalign{\smallskip} \hline \noalign{\smallskip}
\scriptsize SC & \scriptsize Staffing Level & \scriptsize Sq-root Coef. & \scriptsize Resource Req. & \scriptsize Cont. Correction\\ \hline \noalign{\smallskip}
\scriptsize0 & \scriptsize$s$  & \scriptsize$c=\frac{s-R}{\sqrt{R}} $  & \scriptsize$R=\frac{\lambda}{\mu}$ & \scriptsize$\Delta  =\frac{0.5}{\sqrt{R}}$ \\ \noalign{\smallskip} \hline \noalign{\smallskip}
\scriptsize1&\scriptsize$s_1=\frac{s\mu}{\theta_1}$& \scriptsize$c_1 = \frac{s_1-R_1}{\sqrt{R_1}}\left(=\sqrt{\frac{\mu}{\theta_1}} c\right)$&\scriptsize$R_1=\frac{\lambda}{\theta_1}$&\scriptsize$\Delta_1 =\frac{0.5}{\sqrt{R_1}}$\\ \noalign{\smallskip}
 & \scriptsize$s_{1+}=s_1+n_1$&\scriptsize$c_{1+} = \frac{s_{1+}-R_1}{\sqrt{R_1}}$&\scriptsize$\left(=\frac{\mu}{\theta_1}R\right)$&\scriptsize$\left(=\sqrt{\frac{\theta_1}{\mu}}\Delta\right)$\\ \noalign{\smallskip} \hline \noalign{\smallskip}
\scriptsize2 &\scriptsize$s_2=\frac{s\mu+n_1\theta_1}{\theta_2}$&\scriptsize$c_2 = \frac{s_2-R_2}{\sqrt{R_2}}\left(=\sqrt{\frac{\theta_1}{\theta_2}} c_{1+}\right)$&\scriptsize$R_2=\frac{\lambda}{\theta_2}$&\scriptsize$\Delta_2 =\frac{0.5}{\sqrt{R_2}}$\\
 & \scriptsize$s_{2+}=s_2+n_2$&\scriptsize$c_{2+} = \frac{s_{2+}-R_2}{\sqrt{R_2}}$&\scriptsize$\left(=\frac{\theta_1}{\theta_2} R_1=\frac{\mu}{\theta_2}R\right)$&\scriptsize$\left(=\sqrt{\frac{\theta_2}{\theta_1}}\Delta_1=\sqrt{\frac{\theta_2}{\mu}}\Delta \right)$\\ \noalign{\smallskip} \hline 
\end{tabular}\\
\label{tab:parameters}
\end{center}
\normalsize
\end{table}

Corresponding to the three subchains (SC) 0 (M/M/s/s queue), 1 (1\textsuperscript{st} stage queue), and 2 (2\textsuperscript{nd} stage queue), we define parameters necessary for Poisson and standard normal representations in Table~\ref{tab:parameters}. Note that the staffing levels $s$, $s_1$, and $s_2$ must be non-negative integers when they appear in a Poisson representation; if they are non-integer, we round them to their nearest integer values. (In other words, a Poisson representation is exact when staffing levels are integers and otherwise is an approximation.) This integer condition can be dropped when we convert to the normal; a major benefit of using the normal representation. In this conversion procedure, we take into consideration the impact of the discrete state space using the continuity correction terms ($\Delta$, $\Delta_1$, $\Delta_2$). These correction terms can be dropped only at the limit of large resource requirements ($R, R_1, R_2 \to \infty$, respectively).

Using Lemma~\ref{Poisson-Normal} and parameters defined in Table~\ref{tab:parameters}, we obtain Proposition~\ref{block-normal}.
\begin{proposition}
\label{block-normal}
The steady-state probabilities of the subchains are approximately expressed in Poisson and normal representations as follows:
\small
\begin{enumerate}
\item M/M/s/s (subchain 0): 
\begin{equation}
\label{pi0}
\frac{1}{\pi _{s}^{0} } =\frac{F_P(s;R)}{f_P(s;R)} \approx \frac{\sqrt {R}}{h(-c-\Delta)} =: \widetilde{h}.
\end{equation}
\item 1\textsuperscript{st} stage queue (subchain 1):
\begin{equation}
r_1 = \frac{\pi^1_{s+n_1}}{\pi^1_{s}} = \frac{f_P(s_{1+};R_1)}{f_P(s_{1};R_1)} \approx \frac{\phi(c_{1+}+\Delta_1)}{\phi(c_{1}+\Delta_1)} =: \widetilde{r}_1,\nonumber \\
\label{pi1}
\frac{1}{\pi _{s}^{1} } - 1 =\frac{1-F_P(s_1;R_1)}{f_P(s_1;R_1)}-\frac{1-F_P(s_{1+};R_1)}{f_P(s_1;R_1)} \approx \sqrt {R_1}\left(\frac{1}{h(c_1+\Delta_1)}-\frac{\widetilde{r}_1}{h(c_{1+}+\Delta_1)}\right)=: \widetilde{h}_1.
\end{equation}
\item 2\textsuperscript{nd} stage queue (subchain 2):
\begin{equation}
r_2 = \frac{\pi^2_{s+n_1+n_2}}{\pi^2_{s+n_1}} = \frac{f_P(s_{2+};R_2)}{f_P(s_{2};R_2)} \approx \frac{\phi(c_{2+}+\Delta_2)}{\phi(c_{2}+\Delta_2)} =: \widetilde{r}_2,\nonumber \\
\label{pi2}
\frac{1}{\pi _{s+n_1}^{2} } - 1 = \frac{1-F_P(s_2;R_2)}{f_P(s_2;R_2)}-\frac{1-F_P(s_{2+};R_2)}{f_P(s_2;R_2)} \approx \sqrt {R_2}\left(\frac{1}{h(c_2+\Delta_2)}-\frac{\widetilde{r}_2}{h(c_{2+}+\Delta_2)}\right)=: \widetilde{h}_2.
\end{equation}
\end{enumerate}
\normalsize
\end{proposition}

\subsection{Performance Measures in Normal Representation}
\label{sec:normal}
Combining Propositions~\ref{str-rep} and \ref{block-normal} and using parameters defined in Table~\ref{tab:parameters}, we derive the approximate representations of our performance measures.
\begin{proposition}
\label{norm-rep}
Performance measures of the two-stage reneging queue are approximately expressed in the standard normal representation as follows:
\small
\begin{gather*}
\frac{1}{\pi_s} = \widetilde{h}+ \widetilde{h}_1+ \widetilde{r}_1 \widetilde{h}_2, \mbox{ }
\frac{P_Q}{\pi_s} = \frac{1}{\pi _s^q}=1+\widetilde{h}_1 + \widetilde{r}_1 \widetilde{h}_2, \mbox{ }
\frac{P_A}{\pi_s} = p\left(\widetilde{h}_1+\widetilde{r}_1 \widetilde{h}_2\right)+1,\\
\frac{L}{\pi_s} = R_1 (p \widetilde{h}_1 + 1-\widetilde{r}_1)+\widetilde{r}_1 R_2 \left[\left(p+\tfrac{n_1}{R_2}-\tfrac{n_1}{R_1}\right) \widetilde{h}_2 + 1-\widetilde{r}_2 \right].
\end{gather*}
\end{proposition}
Proposition~\ref{norm-rep} can represent both finite and infinite queues as well as the Erlang A model by choosing parameters appropriately. Notice that for $i\in \{1,2\}$, ${n_i \to \infty} \iff {\widetilde{r}_i=0}$; and ${n_i=0} \iff {\widetilde{r}_i=1, \widetilde{h}_i=0}$. The above formulas reduce to the Erlang A formulas either by setting ${\theta_1=\theta_2}$ and ${\widetilde{r}_2=0}$ or by setting ${\widetilde{r}_1=1, \widetilde{h}_1=0}$, and ${\widetilde{r}_2=0}$. Proposition~\ref{norm-rep} can be extended to more stages if desired (see Appendix \ref{three-stage}). Thanks to the linearity of Proposition~\ref{tpt}, Proposition~\ref{norm-rep} shows how each subchain (i.e., $\widetilde{h}$, $\widetilde{h}_1$, $\widetilde{h}_2$) contributes to the performance measures.

Proposition~\ref{norm-rep} can be used to derive asymptotic formulas; such formulas can be useful, but they obscure the impact of the finite 1\textsuperscript{st} stage queue and thus diminish the value of our two-stage model. However, the contrast between asymptotic formulas and our results illustrates some key elements of our results. For example, taking the Halfin-Whitt scaling (i.e., taking a limit of large $R$ and $s$ while maintaining $s=R+c\sqrt{R}$ with fixed $c$) would cause the constant term 1 to vanish from $\frac{P_Q}{\pi_s}$ (i.e., drop ${\pi_s}$ from $P_Q$) as well as all $\Delta$ terms (due to large $R_1$ and $R_2$), leading to a well-known asymptotic formula for the Erlang A model with a single reneging rate $\theta$ in \cite{garnett2002designing}: $P_Q \approx \left[1+\frac{h(c \sqrt{\frac{\mu}{\theta}})}{\sqrt{\frac{\mu}{\theta}} h(-c)}\right]^{-1}$. The difference from $P_Q$ in Proposition~\ref{norm-rep} is the $\pi_s$ and $\Delta$ terms, the correction necessary to accurately represent the discrete state space of MC systems, both of which become negligible in the limit for systems with a large resource requirement, but are critical for non-asymptotic systems.

Propositions \ref{block-normal} and \ref{norm-rep} together with Table~\ref{tab:parameters} provide important insights into our model. Quantitatively speaking, a dominant subchain exhibits a larger parameter, either $\widetilde{h}$, $\widetilde{h}_1$, or $\widetilde{h}_2$, through which it has the correspondingly largest impact on the performance measures, following Proposition~\ref{norm-rep}; but how can we identify which subchain is dominant? The key observation is that, in our approximation scheme, the number of customers in the system is distributed following a normal distribution, which is extended over three subchains. This distribution changes its z-score range in each subchain, but its sign does not change at the boundary of two subchains: subchain 0 ranges up to $c$, subchain 1 ranges from $c_1(=\sqrt{\frac{\mu}{\theta_1}} c)$ to $c_{1+}$, and subchain 2 ranges from $c_2 (=\sqrt{\frac{\theta_1}{\theta_2}} c_{1+})$ to $c_{2+}$ (here, we ignore all continuity correction terms to simplify the argument). Thus, the highest peak of the distribution (corresponding to $z=0$) must fall into one of the subchains (or possibly none of them). We can reasonably assume that if one subchain contains both the highest peak of the distribution and a large portion of distribution, as is common, then this subchain should play an important role in the performance measures. Using this insight, we can make rough capacity planning decisions without evaluating the performance measures. We discuss this more in \S\ref{sec:performance}.

\section{Numerical Experiments}
\label{sec:numerical}
We first demonstrate the accuracy of our approximate formulas and then examine the impact of selected parameters of the two-stage reneging model.
\begin{figure}
    \caption{Comparison between our approximation and exact values for different 1\textsuperscript{st} stage reneging rates $\theta_1$.}
    \label{fig:theta1}
    \centering
    \begin{subfigure}[]{0.32\textwidth}
     \centering
        \includegraphics[scale=0.5]{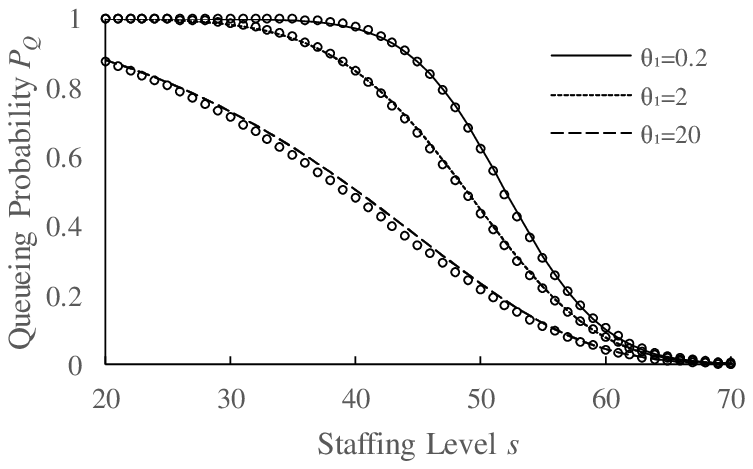}
        \caption{Queueing Probability $P_Q$}
        \label{fig:theta1PQ}
   \end{subfigure}
   \begin{subfigure}[]{0.335\textwidth}
     \centering
        \includegraphics[scale=0.5]{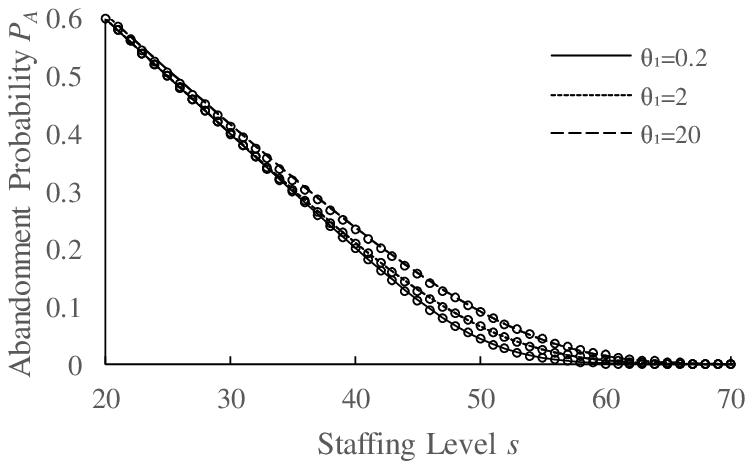}
        \caption{Abandonment Probability $P_A$}
        \label{fig:theta1PA}
    \end{subfigure}
     \begin{subfigure}[]{0.32\textwidth}
     \centering
        \includegraphics[scale=0.5]{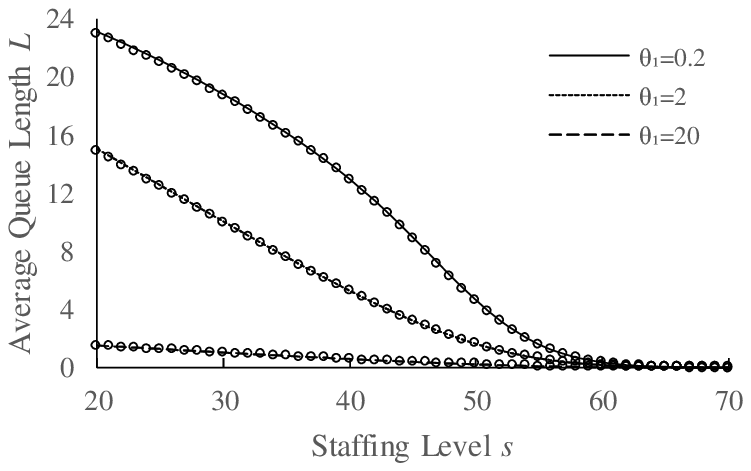}
        \caption{Average Queue Length $L$}
        \label{fig:theta1L}
    \end{subfigure}
\begin{flushleft}
\hspace{-0.2cm} \footnotesize \emph{Notes}: $\lambda=50, \mu=1, n_1=10, n_2=20, \theta_2=2$. $s$ and $\theta_1$ are variables. The open circles in the figures correspond to exact values.
\end{flushleft}
\end{figure}

\subsection{Comparison with Exact Values}
\label{sec:comparison}
\begin{figure}
    \caption{Comparison between our approximation and exact values for different 2\textsuperscript{st} stage reneging rates $\theta_2$.}
    \label{fig:theta2}
    \centering
    \begin{subfigure}[]{0.32\textwidth}
     \centering
        \includegraphics[scale=0.5]{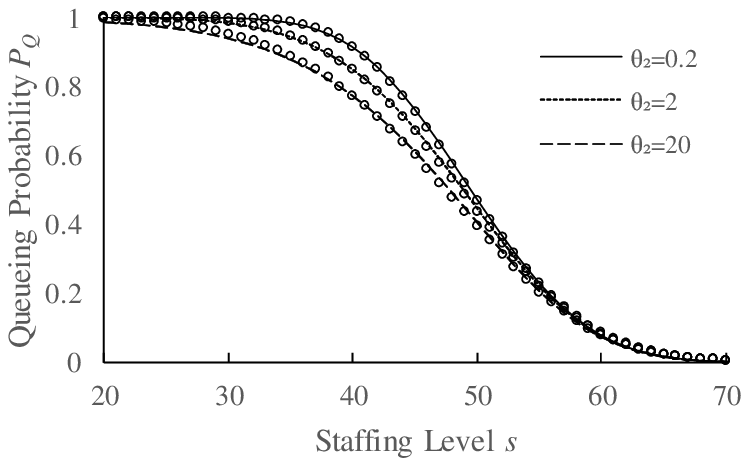}
        \caption{Queueing Probability $P_Q$}
        \label{fig:theta2PQ}
   \end{subfigure}
   \begin{subfigure}[]{0.335\textwidth}
     \centering
        \includegraphics[scale=0.5]{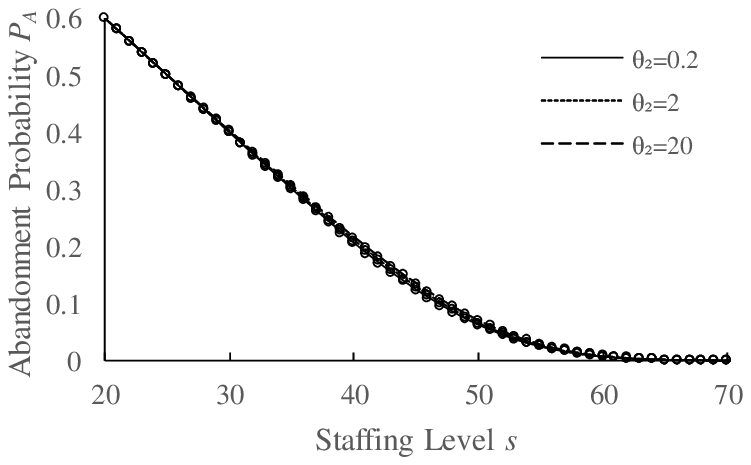}
        \caption{Abandonment Probability $P_A$}
        \label{fig:theta2PA}
    \end{subfigure}
     \begin{subfigure}[]{0.32\textwidth}
     \centering
        \includegraphics[scale=0.5]{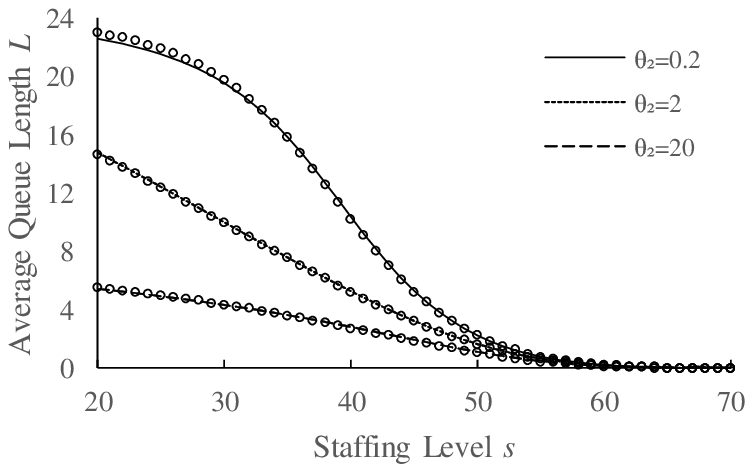}
        \caption{Average Queue Length $L$}
        \label{fig:theta2L}
    \end{subfigure}
\begin{flushleft}
\hspace{-0.2cm} \footnotesize \emph{Notes}: $\lambda=50, \mu=1, n_1=5, n_2=20, \theta_1=2$. $s$ and $\theta_2$ are variables; the lines for $P_A$ are hardly distinguishable due to the the insensitivity of $P_A$ to $\theta_2$ in this example. The open circles in the figures correspond to exact values.
\end{flushleft}
\end{figure}

Fig.~\ref{fig:theta1} compares the exact values with our approximation for fixed $\lambda, \mu, n_1, n_2$, and $\theta_2$. (We tested many parameter settings; results were similar.) Fig.~\ref{fig:theta1} shows staffing level $s$ on the horizontal axis and different lines for different 1\textsuperscript{st} stage reneging rates $\theta_1$. Similarly, Fig.~\ref{fig:theta2} compares the exact values with our approximation when staffing level $s$ changes for different 2\textsuperscript{nd} stage reneging rates $\theta_2$. Tables~\ref{tab:error1} and \ref{tab:error2} show the errors of our approximation corresponding to Figs.~\ref{fig:theta1} and \ref{fig:theta2}, respectively; we present both absolute error (denoted as Abs.; defined as $\emph{exact}-\emph{approximation}$) and relative error (denoted as Rel.; defined as ($\emph{exact}-\emph{approximation})/\emph{exact}$).
\begin{table}
\begin{center}
\caption{Absolute error and relative error of our approximation when varying $\theta_1$ and $s$.}
\scriptsize
\begin{tabular}{lrrrrrrrr}
\noalign{\smallskip} \hline
   & $\theta_1$     & $s$   & 20           & 30          & 40          & 50         & 60          & 70   \\ \hline
$P_Q$ & \multicolumn{1}{r}{0.2} & Abs. & 4.25E-06  & 3.83E-04  & 4.05E-03  & -7.86E-03 & 5.83E-03 & 1.76E-03  \\
   &                         & Rel. & 0.00\%    & 0.04\%    & 0.41\%    & -1.26\%   & 5.47\%   & 37.80\%   \\
   & \multicolumn{1}{r}{2}   & Abs. & 4.96E-04  & 3.81E-03  & -1.11E-03 & -1.10E-02 & 3.69E-03 & 1.50E-03  \\
   &                         & Rel. & 0.05\%    & 0.39\%    & -0.13\%   & -2.52\%   & 4.53\%   & 36.92\%   \\
   & \multicolumn{1}{r}{20}  & Abs. & -6.09E-03 & -1.82E-02 & -2.52E-02 & -1.55E-02 & 4.97E-04 & 9.20E-04  \\
   &                         & Rel. & -0.69\%   & -2.54\%   & -5.24\%   & -7.09\%   & 1.10\%   & 34.52\%   \\ \hline
$P_A$ & \multicolumn{1}{r}{0.2} & Abs. & -2.36E-07 & -2.46E-05 & -4.97E-04 & 1.06E-04  & 4.27E-04 & 2.49E-05  \\
   &                         & Rel. & 0.00\%    & -0.01\%   & -0.25\%   & 0.24\%    & 17.32\%  & 49.66\%   \\
   & \multicolumn{1}{r}{2}   & Abs. & -3.27E-05 & -4.09E-04 & -1.16E-03 & 7.15E-05  & 9.54E-04 & 1.29E-04  \\
   &                         & Rel. & -0.01\%   & -0.10\%   & -0.55\%   & 0.11\%    & 11.79\%  & 44.95\%   \\
   & \multicolumn{1}{r}{20}  & Abs. & -2.61E-03 & -3.02E-03 & -2.37E-03 & 1.28E-04  & 1.74E-03 & 3.62E-04  \\
   &                         & Rel. & -0.43\%   & -0.73\%   & -1.01\%   & 0.14\%    & 10.66\%  & 42.75\%   \\ \hline
$L$  & \multicolumn{1}{r}{0.2} & Abs. & -6.54E-02 & -2.54E-02 & 2.73E-02  & -6.07E-02 & 4.55E-02 & 4.96E-03  \\
   &                         & Rel. & -0.28\%   & -0.14\%   & 0.21\%    & -1.32\%   & 11.12\%  & 45.94\%   \\
   & \multicolumn{1}{r}{2}   & Abs. & -1.97E-02 & -1.21E-02 & -2.92E-02 & 1.79E-03  & 2.39E-02 & 3.23E-03  \\
   &                         & Rel. & -0.13\%   & -0.12\%   & -0.56\%   & 0.11\%    & 11.79\%  & 44.95\%   \\
   & \multicolumn{1}{r}{20}  & Abs. & -6.50E-03 & -7.55E-03 & -5.91E-03 & 3.20E-04  & 4.36E-03 & 9.05E-04  \\
   &                         & Rel. & -0.43\%   & -0.73\%   & -1.01\%   & 0.14\%    & 10.66\%  & 42.75\%  \\ \hline
\end{tabular}
\label{tab:error1}
\begin{flushleft}
\hspace{-0.2cm} \footnotesize \emph{Notes}: $\lambda=50, \mu=1, n_1=10, n_2=20, \theta_2=2$. $\theta_1=0.2, 2, 20$ and $s=20$ to 70.
\end{flushleft}
\end{center}
\normalsize
\end{table}
The results demonstrate that our approximation is accurate, typically with small to very small absolute errors for a wide range of parameters: the reneging rate ranging from patient ($\theta_i=0.2$) to impatient ($\theta_i=20$) and the service capacity ranging from insufficient ($s \ll \lambda/\mu=50$) to ample ($s > \lambda/\mu=50$).
\begin{table}
\begin{center}
\caption{Absolute error and relative error of our approximation when varying $\theta_2$ and $s$.}
\scriptsize
\begin{tabular}{lrrrrrrrr}
\noalign{\smallskip} \hline
   & $\theta_2$     & $s$   & 20           & 30          & 40          & 50         & 60          & 70   \\ \hline
$P_Q$ & \multicolumn{1}{r}{0.2} & Abs. & 1.31E-05  & 8.05E-04  & -6.64E-04 & -8.38E-03 & 4.37E-03 & 1.62E-03   \\
   &                         & Rel. & 0.00\%    & 0.08\%    & -0.07\%   & -1.80\%   & 5.22\%   & 39.69\%    \\
   & \multicolumn{1}{r}{2}   & Abs. & 5.04E-04  & 3.80E-03  & -1.13E-03 & -1.10E-02 & 3.69E-03 & 1.50E-03   \\
   &                         & Rel. & 0.05\%    & 0.38\%    & -0.13\%   & -2.52\%   & 4.53\%   & 36.92\%    \\
   & \multicolumn{1}{r}{20}  & Abs. & 6.53E-03  & 1.11E-02  & -2.01E-03 & -1.45E-02 & 2.29E-03 & 1.43E-03   \\
   &                         & Rel. & 0.66\%    & 1.17\%    & -0.26\%   & -3.69\%   & 2.99\%   & 36.13\%    \\ \hline
$P_A$ & \multicolumn{1}{r}{0.2} & Abs. & -7.46E-07 & -6.29E-05 & -6.63E-04 & -1.76E-04 & 7.98E-04 & 7.92E-05   \\
   &                         & Rel. & 0.00\%    & -0.02\%   & -0.32\%   & -0.28\%   & 10.47\%  & 28.82\%    \\
   & \multicolumn{1}{r}{2}   & Abs. & -3.32E-05 & -4.09E-04 & -1.16E-03 & 7.17E-05  & 9.54E-04 & 1.29E-04   \\
   &                         & Rel. & -0.01\%   & -0.10\%   & -0.55\%   & 0.11\%    & 11.79\%  & 44.95\%    \\
   & \multicolumn{1}{r}{20}  & Abs. & -4.78E-04 & -1.32E-03 & -1.71E-03 & 3.96E-04  & 1.28E-03 & 1.58E-04   \\
   &                         & Rel. & -0.08\%   & -0.33\%   & -0.79\%   & 0.56\%    & 13.91\%  & 47.70\%    \\ \hline
$L$  & \multicolumn{1}{r}{0.2} & Abs. & 4.51E-01  & 2.84E-01  & -1.04E-01 & 5.59E-02  & 3.28E-02 & -1.15E-02  \\
   &                         & Rel. & 1.96\%    & 1.43\%    & -1.01\%   & 2.46\%    & 14.43\%  & -152.68\%  \\
   & \multicolumn{1}{r}{2}   & Abs. & -6.81E-02 & -2.91E-02 & -3.09E-02 & 1.72E-03  & 2.39E-02 & 3.23E-03   \\
   &                         & Rel. & -0.46\%   & -0.29\%   & -0.59\%   & 0.10\%    & 11.79\%  & 44.95\%    \\
   & \multicolumn{1}{r}{20}  & Abs. & 9.15E-02  & 2.98E-02  & -5.59E-02 & -4.60E-02 & 9.98E-03 & 2.64E-03   \\
   &                         & Rel. & 1.65\%    & 0.68\%    & -2.02\%   & -4.31\%   & 6.31\%   & 41.79\%   \\   \hline
\end{tabular}
\label{tab:error2}
\begin{flushleft}
\hspace{-0.2cm} \footnotesize \emph{Notes}: $\lambda=50, \mu=1, n_1=5, n_2=20, \theta_1=2$. $\theta_2=0.2, 2, 20$ and $s=20$ to 70.
\end{flushleft}
\end{center}
\normalsize
\end{table}
This small absolute error nevertheless yields a large relative error when capacity is ample ($s > \lambda/\mu$) since \emph{exact} values of performance measures of interest, the denominator of our statistic, approach zero as $s$ increases.
\begin{figure}
    \caption{Impact of queueing capacities $n_1$ and $n_2$ on performance measures.}
    \label{fig:ex2}
    \centering
    \begin{subfigure}[]{0.32\textwidth}
     \centering
        \includegraphics[scale=0.5]{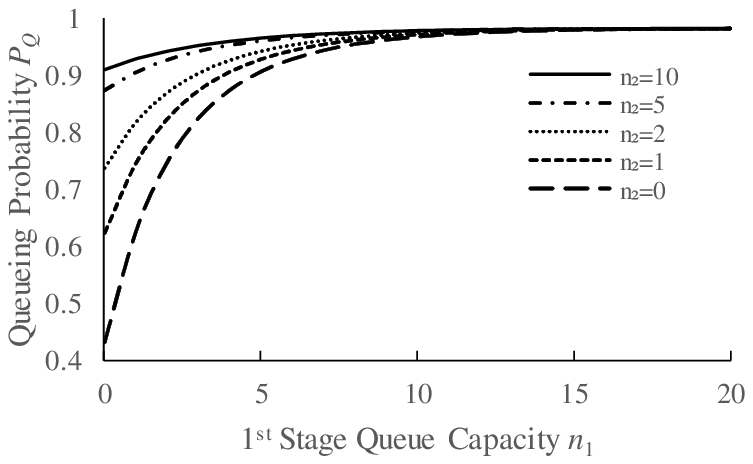}
        \caption{Queueing Probability $P_Q$}
        \label{fig:ex2PQ}
   \end{subfigure}
   \begin{subfigure}[]{0.335\textwidth}
     \centering
        \includegraphics[scale=0.5]{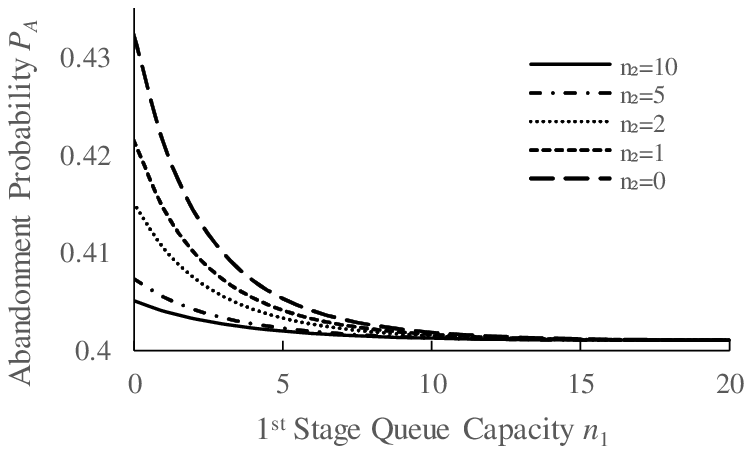}
        \caption{Abandonment Probability $P_A$}
        \label{fig:ex2PA}
    \end{subfigure}
     \begin{subfigure}[]{0.32\textwidth}
     \centering
        \includegraphics[scale=0.5]{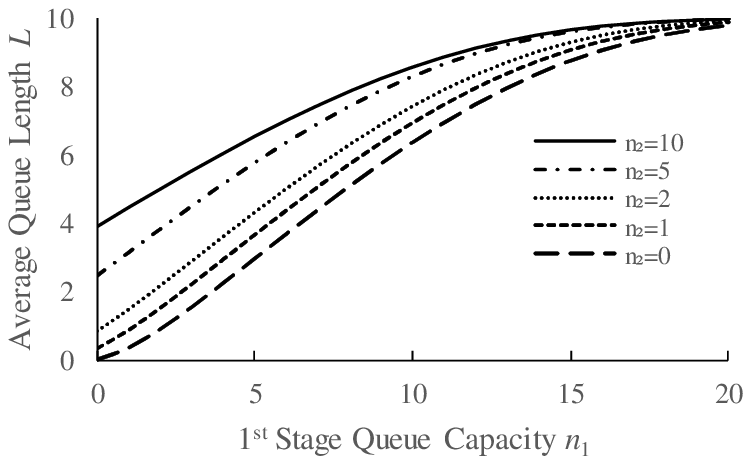}
        \caption{Average Queue Length $L$}
        \label{fig:ex2L}
    \end{subfigure}
\begin{flushleft}
\hspace{1.5cm} \footnotesize \emph{Note}: $\lambda=50, \mu=1, s=30, \theta_1=2, \theta_2=5$. $n_1$ and $n_2$ are variables.
\end{flushleft}
\end{figure}

Finally, we notice from Figs.~\ref{fig:theta1PA} and \ref{fig:theta2PA} that $P_A$ is largely insensitive to reneging rates. This is because customer abandonment is independent from reneging rates in the two extremes: When servers are too slow to handle arrivals (i.e., $s\mu<\lambda$), all customers not served must abandon (at rate $\lambda-s\mu$), and when server capacity is ample, almost nobody abandons. Figs.~\ref{fig:theta1PA} and \ref{fig:theta2PA} confirm this view: $P_A$ approaches $p=(\lambda-s\mu)/\lambda$, a linear function of $s$, when $s$ is small (in the heavy-traffic limit; $s <40$ in our example), while $P_A$ is close to zero when $s$ is large (in the light-traffic limit; $s > 60$ in our example). The reneging rates have some impact on $P_A$ only when $s \approx R$ ($40 \lesssim s \lesssim 60$ in our example).

\subsection{Performance Measures' Dependence on Parameters}
\label{sec:performance}
\begin{figure}
    \caption{Impact of reneging rates $\theta_1$ and $\theta_2$ on performance measures.}
    \label{fig:ex1}
    \centering
    \begin{subfigure}[]{0.32\textwidth}
     \centering
        \includegraphics[scale=0.5]{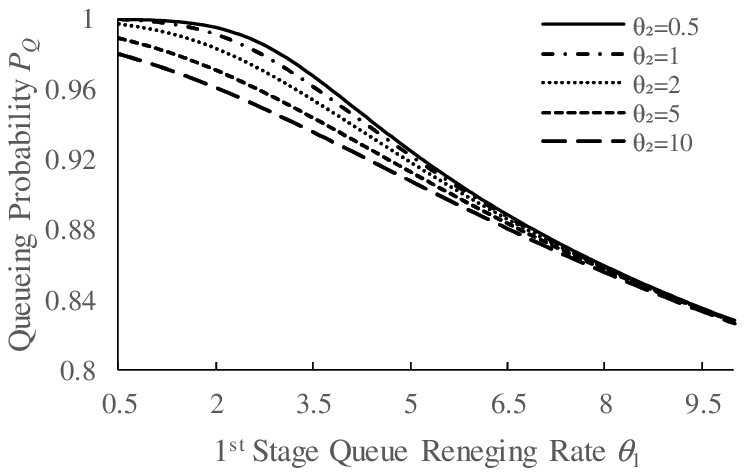}
        \caption{Queueing Probability $P_Q$}
        \label{fig:ex1PQ}
   \end{subfigure}
   \begin{subfigure}[]{0.335\textwidth}
     \centering
        \includegraphics[scale=0.5]{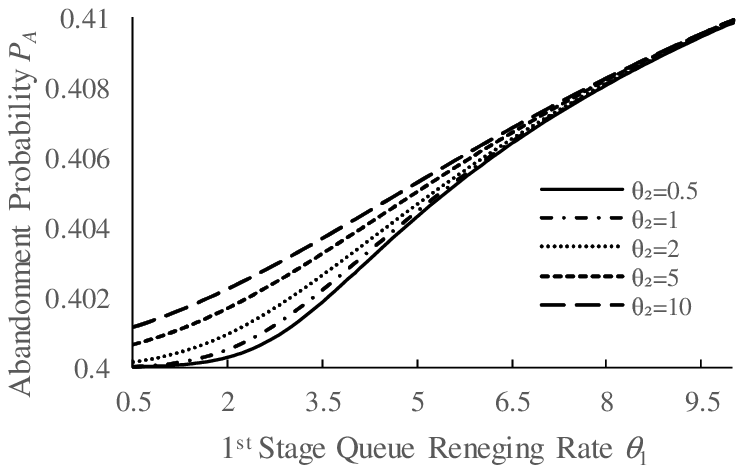}
        \caption{Abandonment Probability $P_A$}
        \label{fig:ex1PA}
    \end{subfigure}
     \begin{subfigure}[]{0.32\textwidth}
     \centering
        \includegraphics[scale=0.5]{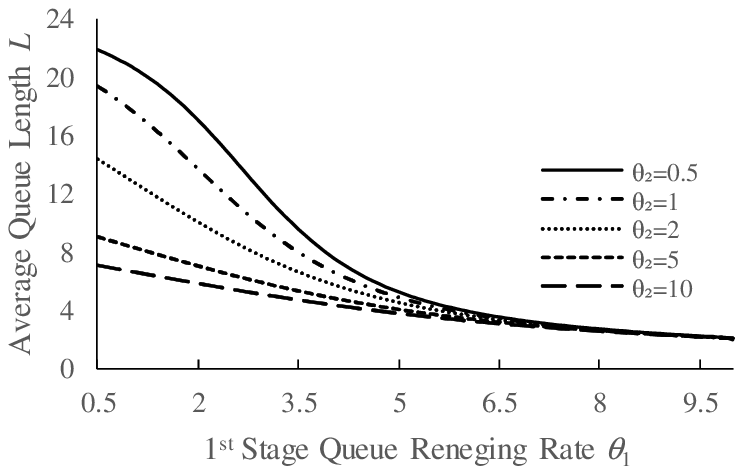}
        \caption{Average Queue Length $L$}
        \label{fig:ex1L}
    \end{subfigure}
\begin{flushleft}
\hspace{1.5cm} \footnotesize \emph{Note}: $\lambda=50, \mu=1, s=30, n_1=6, n_2=20$. $\theta_1$ and $\theta_2$ are variables.
\end{flushleft}
\end{figure}

\begin{figure}
    \caption{Impact of 2\textsuperscript{nd} stage reneging rate $\theta_2$ and 1\textsuperscript{st} stage queue capacity $n_1$ on performance measures.}
    \label{fig:ex3}
    \centering
    \begin{subfigure}[]{0.32\textwidth}
     \centering
        \includegraphics[scale=0.5]{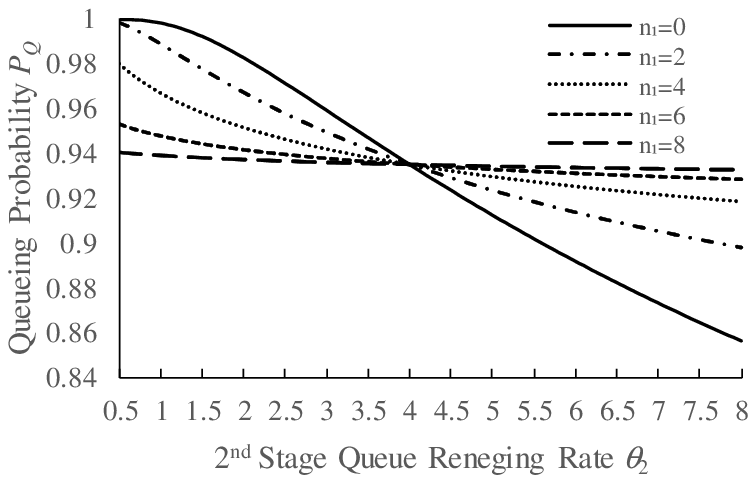}
        \caption{Queueing Probability $P_Q$}
        \label{fig:ex3PQ}
   \end{subfigure}
   \begin{subfigure}[]{0.335\textwidth}
     \centering
        \includegraphics[scale=0.5]{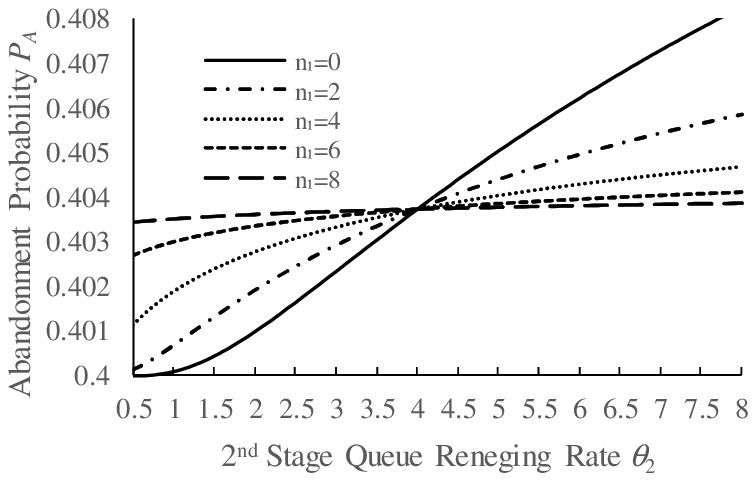}
        \caption{Abandonment Probability $P_A$}
        \label{fig:ex3PA}
    \end{subfigure}
     \begin{subfigure}[]{0.32\textwidth}
     \centering
        \includegraphics[scale=0.5]{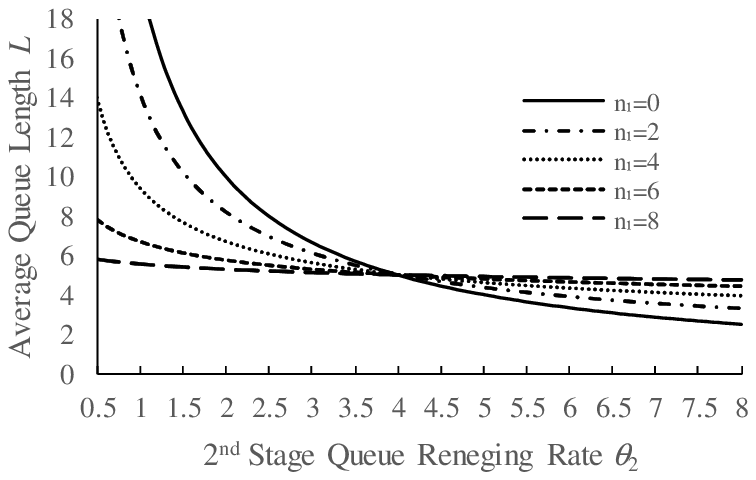}
        \caption{Average Queue Length $L$}
        \label{fig:ex3L}
    \end{subfigure}
\begin{flushleft}
\hspace{1.5cm} \footnotesize \emph{Note}: $\lambda=50, \mu=1, s=30, \theta_1=4, n_2=\infty$. $\theta_2$ and $n_1$ are variables.
\end{flushleft}
\end{figure}

It is intuitive that if customers rarely exceed the 1\textsuperscript{st} stage capacity---for example, due to large $n_1$ and/or large $\theta_1$---then the system will become insensitive to $n_2$ and $\theta_2$. Numerical experiments verify this: Figs.~\ref{fig:ex2} and \ref{fig:ex1} show how $n_1$ and $\theta_1$ (horizontal axis) impact performance measures for different $n_2$ and $\theta_2$, respectively; we observe all lines converge as $n_1$ and $\theta_1$ increase. To illustrate how the 1\textsuperscript{st} stage queue affects the system's dependence on the 2\textsuperscript{nd} stage queue, we conduct the following experiment: Fig.~\ref{fig:ex3} shows the impact of $\theta_2$ (horizontal axis) on performance measures for different $n_1$; a steeper line corresponds to higher impact of $\theta_2$. We observe that, in this example, at least about $n_1=8$ is necessary to make the system insensitive to the change of $\theta_2$. A real-world implication of this result is that, by conducting a similar analysis on a spreadsheet, a restaurant owner could decide how many seats to provide for customers waiting inside (1\textsuperscript{st} stage queue) to make the quality of service robust to the weather, which directly affects the level of patience for customers waiting outside (2\textsuperscript{nd} stage queue).

Before we conclude this section, we provide further insight into this insensitivity. First, notice that the 1\textsuperscript{st} stage queue is (approximately) represented by the standard normal distribution ranging from $c_1(=\frac{s\mu - \lambda}{\sqrt{\lambda \theta_1}})$ to $c_{1+}(=\frac{s\mu+n_1 \theta_1 - \lambda}{\sqrt{\lambda \theta_1}})$. To contain most waiting customers within the 1\textsuperscript{st} stage queue, its distribution should cover a wide range of the standard normal distribution including the center of the distribution. Since $c_1<0$ (i.e., $\lambda > s\mu$) in the case of interest (otherwise, there are only a small number of waiting customers), we require $c_{1+}  \ge z$, where $z$ is a decision variable representing a certain positive threshold on the standard normal distribution. For example, if we use $z=1$, we obtain $n_1 \ge15$ in Fig.~\ref{fig:ex2}; $\theta_1 \ge 6.3$ in Fig.~\ref{fig:ex1}; and $n_1 \ge 8.5$ (or 9 since $n_1$ should be an integer) in Fig.~\ref{fig:ex3}; these estimates turn out to be sufficiently good. Although this rule of thumb cannot replace accurate numerical calculations based on Proposition~\ref{norm-rep}, it allows practitioners to make ballpark parameter estimates for containing most waiting customers in the 1\textsuperscript{st} stage queue.

\section{Conclusion and future work}
\label{sec:conclusion}
In this paper we study a single station two-stage reneging model with finite or infinite queue capacity. Our two-stage model has wide applicability in the real world; for example, practitioners can utilize our formulas to make capacity decisions efficiently without solving Markov chains or conducting simulations.

There are several possible extensions of our results. First, our two-stage model can be straightforwardly extended to three or more stages if further accuracy is needed. Second, our paper assumes a constant arrival rate with no balking, but extension to a stage-dependent arrival rate (i.e., stage-dependent balking) is straightforward. Third, in addition to the three performance measures we analyze in this paper, we can also analyze other performance measures such as waiting time distribution. Our two-stage reneging model is simple, but more versatile than the Erlang A model that has inspired researchers for several decades; we hope our two-stage model opens a new avenue of research on service systems with impatient customers.

\appendixtitleon
\begin{appendices}
\section{Proof of Proposition~\ref{tpt}}
Since the MC we consider is a birth-death type, all decomposed subchains maintain a stationary distribution proportional to the full MC: for any subchain $j$, $\pi_k^j \propto \pi_k, \forall k \in A_j$, or equivalently, $\pi_k/\pi^j_{k}= \pi_{k'}/\pi^j_{k'}, \forall k,k' \in A_j$. Also, note that $A_0 \cap A_q = \{s\}$. Define $\mathbb{E}[{f(X)}]=\sum_{k \in \Omega} f(k)\pi_k$ and $\mathbb{E}_j[{f(X)}]=\sum_{k \in A_j} f(k)\pi^j_k$. Applying the total expectation theorem to the full MC, we derive
\small
\begin{align*}
\mathbb{E}[f(X)]&=\sum_{k \in A_0} f(X) \pi_k + \sum_{k \in A_q} f(X) \pi_k - f(s) \pi_s\\
	&=\sum_{k \in A_0} f(X) \pi^0_k \frac{\pi_s}{\pi^0_s} + \sum_{k \in A_q} f(X) \pi^q_k \frac{\pi_s}{\pi^q_s}- f(s) \pi_s\\
	&=\mathbb{E}_0[f(X)] \frac{\pi_s}{\pi^0_s}  + \mathbb{E}_q[f(X)] \frac{\pi_s}{\pi^q_s} - f(s) \pi_s,
\end{align*}
\normalsize
from which we obtain \eqref{tpt1}. The derivation of \eqref{tpt2} is almost identical. Note that $A_1 \cap A_2 = \{s+n_1\}$. By decomposing subchain $q$ into subchains 1 and 2, we obtain
\small
\begin{equation*}
\frac{\mathbb{E}_q[f(X)]}{\pi_{s+n_1}^q}=\frac{\mathbb{E}_1[f(X)]}{\pi^1_{s+n_1}}  + \frac{\mathbb{E}_2[f(X)]}{\pi^2_{s+n_1}} - f(s+n_1).
\end{equation*}
\normalsize
Multiplying both sides of the equation by ${r_1=\frac{\pi^1_{s+n_1}}{\pi^1_s}}$, we obtain \eqref{tpt2}. Finally, \eqref{tpt3} directly follows from \eqref{tpt1} and \eqref{tpt2}.

\section{Proof of Proposition~\ref{str-rep}}
For \eqref{pi}, we set $f(X)=1$ in Proposition~\ref{tpt} and obtain $\frac{1}{\pi_s}=\frac{1}{\pi^0_s}+\frac{1}{\pi^q_s}-1$ and $\frac{1}{\pi^q_s}=\frac{1}{\pi^1_s}+r_1(\frac{1}{\pi^2_{s+n_1}}-1)$, from which we obtain the result.
For \eqref{pq}, we set $f(X)=\bold{1}_{A_q}$, an indicator function. Note that $P_Q=E[\bold{1}_{A_q}]$, $E_0[\bold{1}_{A_q}]=\pi^0_s$, $E_q[\bold{1}_{A_q}]=1$, and $E_s[\bold{1}_{A_q}]=1$. Thus, we obtain $\frac{P_Q}{\pi_s}=\frac{E[\bold{1}_{A_q}]}{\pi_s}=\frac{1}{\pi^q_s}=\frac{1}{\pi^1_s}+r_1(\frac{1}{\pi^2_{s+n_1}}-1)$.
For \eqref{pa}, we set $f(X)=N_A$, a random variable representing the steady-state number of customer abandonments (via reneging and blocking) per unit time. Note that $P_A=\frac{E[N_A]}{\lambda}$, $E_0[N_A]=\lambda \pi^0_s$, $E_s[N_A]=\lambda$, and $E_q[N_A]=\lambda-s \mu(1-\pi^q_s)$ (the flow conservation law). Thus, we obtain $\frac{P_A}{\pi_s}=\frac{E[N_A]}{\lambda \pi_s}=p(\frac{1}{\pi^q_s}-1)+1$, where $p=1-\frac{s\mu}{\lambda}$.
Finally, as for \eqref{el}, we set ${f(X)=N}$, a random variable representing the steady-state number of customers in queue. Let ${L=E[N]}$, ${L_1=E_1[N]}$, and ${L_2=E_2[N]}$. We obtain ${\frac{L}{\pi_s}=\frac{E[N]}{\pi_s}=\frac{L_1}{\pi^1_s}+r_1(\frac{L_2}{\pi^2_{s+n_1}}-n_1)}$. The remaining task is to find $L_1$ and $L_2$. We use the flow conservation law for subchains 1 and 2: $\lambda (1-\pi^1_{s+n_1})=s\mu(1-\pi^1_s)+\theta_1 L^1$ and $\lambda (1-\pi^2_{s+n_1+n_2})=(s\mu+n_1\theta_1)(1-\pi^2_{s+n_1})+\theta_2 (L^2-n_1)$, from which we obtain the expressions for $L_1$ and $L_2$:
\small
\begin{align*}
\frac{L^1}{\pi^1_s}&=\frac{\lambda}{\theta_1}\left[p\left(\frac{1}{\pi^1_s}-1\right)+1-r_1\right],\\
\frac{L^2}{\pi^2_{s+n_1}}-n_1&=\frac{\lambda}{\theta_2}\left[\left(p-\frac{n_1(\theta_2-\theta_1)}{\lambda}\right)\left(\frac{1}{\pi^2_{s+n_1}}-1\right)+1-r_2\right],
\end{align*}
\normalsize
where ${r_1=\frac{\pi^1_{s+n_1}}{\pi^1_{s}}}$ and ${r_2=\frac{\pi^2_{s+n_1+n_2}}{\pi^2_{s+n_1}}}$.

\section{Proof of Corollary \ref{rules}}
The results directly follow from Proposition~\ref{str-rep}. We obtain \eqref{pa1} by eliminating $\pi^q_s$ in \eqref{pq} and \eqref{pa}. We obtain \eqref{pa3} by setting $\theta=\theta_1=\theta_2$ in \eqref{el} and comparing it with \eqref{pq}.

\section{Proof of Lemma~\ref{Poisson-Normal}}
The Poisson-to-Normal conversion is based on the Central Limit Theorem and is described in many textbooks \cite[e.g.,][]{feller1968probability}. The approximation achieves high accuracy when the continuity correction term $\Delta$ is small ($R$ is large); at the limit of small $\Delta$ (large $R$), the error of the approximation approaches zero. For \eqref{cdf}, we make a discrete-to-continuous conversion from Poisson to normal with a continuity correction ``+0.5" following \cite{tijms1986stochastic}. We then convert normal to standard normal using $c=\frac{s-R}{\sqrt R}$ and $\Delta=\frac{0.5}{\sqrt R}$: ${F_P}(s;R) \approx  \Phi \left(\frac{(s+0.5)-R}{\sqrt{R}}\right)= \Phi \left(c + \Delta\right).$ For \eqref{pdf}, using the result above and the assumption that $\Delta$ is sufficiently small, we obtain ${f_P}(s;R) = {F_P}(s;R) - {F_P}(s-1;R) \approx \Phi (c + \Delta) - \Phi (c - \Delta) \approx \frac{\phi (c + \Delta)}{\sqrt R}.$ For \eqref{hazard1} and \eqref{hazard2}, we use the above results and the definition of the hazard function for the standard normal distribution. We obtain $\frac{f_P(s;R)}{1-F_P(s;R)} \approx \frac{\phi (c+\Delta)}{\sqrt{R} \left(1-\Phi (c+\Delta)\right)} =\frac{h(c+\Delta)}{\sqrt{R}} \text{ and } \frac{f_{P}(s;R)}{F_{P}(s;R)} \approx \frac{\phi (c+\Delta)}{\sqrt{R} \Phi (c+\Delta)} =\frac{h(-c-\Delta)}{\sqrt{R}}.$

\section{Proof of Proposition~\ref{block-normal}}
We discuss each subchain separately. We use parameters defined in Table~\ref{tab:parameters}. We assume that all staffing levels $s$, $s'$, and $s''$ are non-negative integers in a Poisson representation; if they are non-integer, we round them to their nearest integer values. This integer condition is dropped when we convert to the normal.

\begin{enumerate}
\item M/M/s/s (subchain 0): This result is familiar \cite[see, e.g.,][]{harchol2013performance}. Let ${R=\frac{\lambda}{\mu}}$ and consider a random variable $X \sim \text{Pois}(R)$; thus Pr${\{X=i\}=\frac{e^{-R}R^i}{i!}}$, ${\forall i \in \mathbb{Z}_{\ge 0}}$. This subchain is a birth-death MC with the total departure rate $k\mu$ at state $k$. Hence, for all $k=1,2,\dots,s $,
\small
$$\pi _k^0 = \pi _{k - 1}^0\frac{\lambda}{k\mu} = \pi _{k - 1}^0\frac{R}{k} = \ldots = \pi _0^0\frac{R^k}{k!} \propto \Pr\{X=k\}.$$
\normalsize
Note that $\pi_k^0 \propto \Pr\{X=k\}$ means the distributions of $\pi _k^0$ and ${\Pr\{X=k\}}$ for $k \in A_0$ are the same except for the normalization constant.
By summing up the terms with respect to $k$ and applying the normalization condition, we obtain
\footnotesize
\begin{equation*}
\frac{1}{\pi_s^0} = \frac{\Pr \{0 \le X \le s\}}{\Pr \{X =s\}} = \frac{F_{P}(s;R)}{f_{P}(s;R)}.
\end{equation*}
\normalsize
Let $c=\frac{s-R}{\sqrt{R}}$ and $\Delta=\frac{0.5}{\sqrt{R}}$. Using \eqref{hazard2}, we obtain \eqref{pi0}.

\item 1\textsuperscript{st} stage queue (subchain 1): Let $R_1=\frac{\lambda}{\theta_1}$, $s_1=\frac{s\mu}{\theta_1}$, and ${s_{1+}=s_1+n_1}$. Consider a random variable $X_1 \sim \text{Pois}(R_1)$; thus Pr${\{X_1=i\}=\frac{e^{-R_1}(R_1)^i}{i!}}$, ${\forall i \in \mathbb{Z}_{\ge 0}}$. This subchain is a birth-death MC with the total departure rate $s\mu+k\theta_1=(\frac{s\mu}{\theta_1}+k)\theta_1=(s_1+k)\theta_1$ at state $s+k$. For all $k=1,2,\dots,n_1$,
\small
\begin{align*}
\pi _{s+k}^1 &= \pi _{s+k-1}^1\frac{\lambda}{(s_1+k)\theta_1} = \pi_{s+k-1}^1\frac{R_1}{s_1+k} = \ldots = \pi _s^1\frac{(R_1)^{s_1+k}}{(s_1+k)!} \frac{s_1!}{(R_1)^{s_1}}\\
&\propto \Pr\{X_1=s_1+k\}.
\end{align*}
\normalsize
By summing up the terms with respect to $k$ and applying the normalization condition, we obtain
\footnotesize
\begin{align*}
\frac{1}{\pi_s^1} &= \frac{\Pr \{s_1 \le X_1 \le s_1+n_1\}}{\Pr \{X_1 =s_1\}}=\frac{\Pr \{X_1 = s_1\}-\Pr \{X_1 \le s_1\}+\Pr \{X_1 \le s_1+n_1\}}{\Pr \{X_1 =s_1\}}\\
&= 1+ \frac{1-F_P(s_1;R_1)}{f_P(s_1;R_1)}-\frac{1-F_P(s_{1+};R_1)}{f_P(s_1;R_1)}\\
&=1+ \frac{1-F_P(s_1;R_1)}{f_P(s_1;R_1)}-\frac{f_P(s_{1+};R_1)}{f_P(s_{1};R_1)}\frac{1-F_P(s_{1+};R_1)}{f_P(s_{1+};R_1)}.
\end{align*}
\normalsize
Let $c_1=\frac{s_1-R_1}{\sqrt{R_1}}$, $c_{1+}=\frac{s_{1+}-R_1}{\sqrt{R_1}}$, and $\Delta_1=\frac{0.5}{\sqrt{R_1}}$. Using \eqref{pdf} and \eqref{hazard1}, we obtain \eqref{pi1}.

\item 2\textsuperscript{nd} stage queue (subchain 2): The argument is almost identical to the case of the 1\textsuperscript{st} stage queue. Let $R_2=\frac{\lambda}{\theta_2}$, ${s_2=\frac{s\mu+n_1\theta_1}{\theta_2}}$, and ${s_{2+}=s_2+n_2}$. Consider a random variable $X_2 \sim \text{Pois}(R_2)$; thus Pr${\{X_2=i\}=\frac{e^{-R_2}(R_2)^i}{i!}}$, ${\forall i \in \mathbb{Z}_{\ge 0}}$. This subchain is a birth-death MC with the total departure rate $s\mu+n_1\theta_1+k\theta_2=(s_2+k)\theta_2$ at state $s+n_1+k$. For all $k=1,2,\dots,n_2$,
\footnotesize
\begin{align*}
\pi _{s+n_1+k}^2 &= \pi _{s+n_1+k-1}^2 \frac{\lambda}{(s_2+k)\theta_2} = \pi_{s+n_1+k-1}^2 \frac{R_2}{s_2+k} = \dots \\
&= \pi _{s+n_1}^2 \frac{(R_2)^{s_2+k}}{(s_2+k)!} \frac{s_2!}{(R_2)^{s_2}}\propto \Pr\{X_2=s_2+k\},
\end{align*}
\normalsize
from which we obtain
\smallskip
\footnotesize
\begin{equation*}
\frac{1}{\pi_{s+n_1}^2} = \frac{\Pr \{s_2 \le X_2 \le s_2+n_2\}}{\Pr \{X_2 =s_2\}}=1+ \frac{1-F_P(s_2;R_2)}{f_P(s_2;R_2)}-\frac{f_P(s_{2+};R_2)}{f_P(s_{2};R_2)}\frac{1-F_P(s_{2+};R_2)}{f_P(s_{2+};R_2)}.
\end{equation*}
\normalsize
\smallskip
Let $c_2=\frac{s_2-R_2}{\sqrt{R_2}}$, $c_{2+}=\frac{s_{2+}-R_2}{\sqrt{R_2}}$, and $\Delta_2=\frac{0.5}{\sqrt{R_2}}$. Using \eqref{pdf} and \eqref{hazard1}, we obtain \eqref{pi2}.
\end{enumerate}

\section{Proof of Proposition~\ref{norm-rep}}
The results directly follow from Propositions~\ref{str-rep} and \ref{block-normal}.

\section{Extension to a Three-Stage Reneging Queue}
\label{three-stage}
To extend the two-stage reneging model to three-stage, we introduce the 3\textsuperscript{rd} subchain with capacity $n_3$ and reneging rate $\theta_3$. Define parameters accordingly: ${R_3=\frac{\lambda}{\theta_3}}$, ${s_3=\frac{s\mu+n_1\theta_1+n_2\theta_2}{\theta_3}}$, ${s_{3+}=s_3+n_3}$, ${c_3=\frac{s_3-R_3}{\sqrt{R_3}}}$, ${c_{3+}=\frac{s_{3+}-R_3}{\sqrt{R_3}}}$, ${\Delta_3=\frac{0.5}{\sqrt{R_3}}}$. Also, we introduce $r_3 = \frac{\pi^3_{s+n_1+n_2+n_3}}{\pi^3_{s+n_1+n_2}} \approx \frac{\phi(c_{3+}+\Delta_3)}{\phi(c_{3}+\Delta_3)} =: \widetilde{r}_3$ and $\frac{1}{\pi _{s+n_1+n_2}^{3}} - 1 \approx \sqrt {R_3}\left(\frac{1}{h(c_3+\Delta_3)}-\frac{r_3}{h(c_{3+}+\Delta_3)}\right)=: \widetilde{h}_3$. It is straightforward to obtain the expressions for the three-stage reneging queue.
\footnotesize
\begin{equation*}
\frac{1}{\pi_s} = \widetilde{h}+ \widetilde{h}_1+ \widetilde{r}_1 \widetilde{h}_2+\widetilde{r}_1 \widetilde{r}_2 \widetilde{h}_3, \mbox{ }
\frac{P_Q}{\pi_s} = 1+\widetilde{h}_1 + \widetilde{r}_1 \widetilde{h}_2+\widetilde{r}_1 \widetilde{r}_2 \widetilde{h}_3, \mbox{ }
\frac{P_A}{\pi_s} = p\left(\widetilde{h}_1+\widetilde{r}_1 \widetilde{h}_2+\widetilde{r}_1 \widetilde{r}_2 \widetilde{h}_3\right)+1,
\end{equation*}
\begin{align*}
\frac{L}{\pi_s} &= R_1 (p \widetilde{h}_1 + 1-\widetilde{r}_1)+\widetilde{r}_1 R_2 \left[(p+\tfrac{n_1}{R_2}-\tfrac{n_1}{R_1}) \widetilde{h}_2 + 1-\widetilde{r}_2 \right]\\
&+\widetilde{r}_1 \widetilde{r}_2 R_3 \left[(p+\tfrac{n_1+n_2}{R_3}-\tfrac{n_1}{R_1}-\tfrac{n_2}{R_2}) \widetilde{h}_3 + 1-\widetilde{r}_3 \right].
\end{align*}
\normalsize
Observe that the formulas we obtain here are simply the formulas in Proposition~\ref{norm-rep} with extra terms added to account for the 3\textsuperscript{rd} subchain. If it is necessary, an extension to four or more stages is also straightforward.
\end{appendices}

\bibliography{abandonment-ref2} 
\end{document}